\documentclass[12pt]{article}
\usepackage{amssymb}
\usepackage{amsmath}
\usepackage{graphics}
\usepackage{natbib}
\usepackage{color}
\bibpunct{(}{)}{;}{x}{,}{,}
\usepackage[colorlinks=true,citecolor=blue]{hyperref}
\newcommand{\ip}[2]{\left <{#1},{#2}\right >}

\newcommand{\bc}{\begin{center}}
\newcommand{\ec}{\end{center}}
\newcommand{\be}{\begin{equation}}
\newcommand{\ee}{\end{equation}}
\newcommand{\bea}{\begin{eqnarray}}
\newcommand{\eea}{\end{eqnarray}}
\newcommand{\bes}{\begin{eqnarray*}}
\newcommand{\ees}{\end{eqnarray*}}
\newcommand{\bed}{\begin{description}}
\newcommand{\eed}{\end{description}}
\newcommand{\bee}{\begin{enumerate}}
\newcommand{\eee}{\end{enumerate}}

\newcommand{\cov}{\mbox{cov}}
\newcommand{{\mbt}}{{\mbox{\boldmath $\beta$}} }
\newcommand{{\mbx}}{{\mbox{\boldmath $\xi$}} }
\newcommand{{\mbp}}{{\mbox{\boldmath $\psi$}} }
\newcommand{{\wmbt}}{\widehat{\mbt}} 
\newcommand{{\wtb}}{\widetilde{\mbt}} 
\newcommand{{\mbv}}{{\mbox{\boldmath $\vartheta$}} }

\newcommand{{\mba}}{{\mbox{\boldmath $\alpha$}} }
\newcommand{{\mbl}}{{\mbox{\boldmath $\lambda$}} }
\newcommand{{\mbL}}{{\mbox{\boldmath $\Lambda$}} }
\newcommand{{\mbO}}{{\mbox{\boldmath $\Omega$}} }

\newcommand{\bS}{\mathbf{\Sigma}}

\newcommand{\bF}{{\bf F}}
\newcommand{\bY}{{\bf Y}}
\newcommand{\bA}{{\bf A}}
\newcommand{\bb}{{\bf b}}
\newcommand{\bR}{{\bf R}}
\newcommand{\bC}{{\bf C}}
\newcommand{\bW}{{\bf W}}
\newcommand{\bV}{{\bf V}}
\newcommand{\bD}{{\bf D}}
\newcommand{\bd}{{\bf d}}
\newcommand{\bJ}{{\bf J}}

\newcommand{\wC}{\widehat{\bC}}
\newcommand{\mB}{\mathcal{B}}

\newcommand{\bI}{{\bf I}}
\newcommand{\bM}{{\bf M}}

\newcommand{\bZ}{{\bf Z}}
\newcommand{\bX}{{\bf X}}
\newcommand{\bg}{{\bf g}}
\newcommand{\bh}{{\bf h}}

\newcounter{assumption}

\newcounter{theorem}

\newcounter{corollary}

\newcounter{lemma}

\newcounter{remark}

\begin{document}

\bc
{\large \bf On Large-Sample Estimation and Testing via Quadratic
Inference Functions for Correlated Data}  
\ec

\bc
{\bf Ramani S.\ Pilla} and {\bf Catherine Loader} \\
{\em Department of Statistics, Case Western Reserve University,
Cleveland, OH 44106} \\
{\em pilla@case.edu \hspace{0.1in} catherine@case.edu} \\ 
December 2005
\ec

\centerline{\bf Abstract}

\cite{hansen:82} proposed a class of {\em generalized method of moments} (GMMs)
for estimating a vector of regression parameters from a set of score
functions. Hansen established that, under certain regularity
conditions, the estimator based on the GMMs is consistent,
asymptotically normal and asymptotically efficient. In the generalized
estimating equation framework, extending the principle of the GMMs to
implicitly estimate the underlying correlation structure leads to a
{\em quadratic inference function} (QIF) for the analysis of
correlated data. The main objectives of this research are to (1)
formulate an appropriate estimated covariance matrix for the set of
extended score functions defining the inference functions; (2) develop
a unified large-sample theoretical framework for the QIF; (3) derive a
generalization of the QIF test statistic for a general linear
hypothesis problem involving correlated data while establishing the
asymptotic distribution of the test statistic under the null and local
alternative hypotheses; (4) propose an iteratively reweighted
generalized least squares algorithm for inference in the QIF
framework; and (5) investigate the effect of basis matrices, defining
the set of extended score functions, on the size and power of the QIF
test through Monte Carlo simulated experiments.

\vspace{0.1in}
\noindent{\em Key Words:} Covariance structure; Extended score;
Generalized estimating equations; Generalized least squares;
Generalized method of moments; Longitudinal data; Quasi-likelihood. 

\section{Introduction}
\label{sec-intro}

Correlated data arises when a response is measured at repeated
instances on a set of subjects within a study design. A canonical
problem is to determine a regression relationship between the measured
responses and a set of covariates. Responses on different subjects are
assumed to be independent, while the repeated measurements on
individual subjects are correlated with an unknown correlation
structure. Any inferential procedure must take account of this
correlation \citep{crowder:90,diggle:94}.

Let $Y_{it}$ be a response, with a corresponding $q$-dimensional
vector of covariates \(\bX_{it}\), measured at the $t$th ($t = 1,
\ldots, n_{i}$) time point on the $i$th ($i = 1, \ldots, N$)
subject. Assuming a generalized linear model for \(Y_{it}\) and
\(\bX_{it}\) yields
\bea
  \label{eq:link}
  E(Y_{it}) = h\left(\bX_{it}^T \, \mbt \right)
\eea
and
\bea
  \label{eq:var} 
\hbox{var}(Y_{it}) = \phi \; v\left(\bX_{it}^T \, \mbt \right),
\eea
where \(\mbt \in \mB\) is a $q$-dimensional vector of unknown
regression parameters, $\phi$ is a dispersion parameter,
\(h(\,\cdot\,)\) is a known inverse link function and $v(\,\cdot\,)$
is a known variance function. The principal objective is to derive
inferential theory for the unknown parameter vector \(\mbt\).

For simplicity of exposition, we assume that each subject is observed
at a common set of time points \(t = 1, \ldots, n\). Let \(\bY_i =
\left(Y_{i1}, \ldots, Y_{in} \right)^{T}\) be the response vector and
$\bh_i = E(\bY_i) = \left\{ h(\bX_{i1}^{T} \, \mbt), \ldots,
h(\bX_{in}^{T} \, \mbt) \right\}^{T}$ be the vector of
means. Furthermore, let the operator \(\nabla\) denote a partial
derivative with respect to the elements of \(\mbt\) so that \(\nabla
\bh_i\) represents the $(n \times q)$ matrix \((\partial
\bh_i/\partial \mbt_1, \ldots, \partial \bh_i/\partial \mbt_q)\) for
each $i = 1, \ldots, N$.  


\subsection{Background}

The quasi-likelihood estimating equation \citep{Wedder:74} for $\mbt$,
under the generalized linear model framework is defined as 
\begin{equation}
  \label{eq:quasi} S(\mbt) := \sum_{i = 1}^{N} \left\{ \nabla
  \bh_{i}^{T} \; \bW_{i}^{-1} \; \left({\bY_{i}} - \bh_{i} \right)
  \right\} = 0,
\end{equation}
where $\bW_i = \cov(\bY_i)$ is an $(n \times n)$ diagonal matrix with
elements determined by the variance function (\ref{eq:var}). When the
independence assumption within a subject for the responses is relaxed,
the matrices \(\bW_i\) are no longer diagonal, instead have an unknown
correlation structure that needs to be incorporated into the model.
In a seminal article, \cite{liang:86} proposed {\em generalized
estimating equations} (GEEs), based on the ingenious idea of using a 
\emph{working correlation matrix} with a nuisance parameter vector to 
simplify $\bW_i$. In particular, they proposed the GEEs based on 
\[
  \bW_i = \bA_i^{1/2} \, \bR(\mba) \, \bA_i^{1/2} \quad \mbox{for}
  \quad i = 1, \ldots, N,
\]
where $\bA_i$ is the diagonal matrix of marginal variance of $\bY_i$
and $\bR(\mba)$ is the {\em working correlation matrix} with an
unknown nuisance parameter vector $\mba$. Specific choices of
\(\bR(\mba)\) correspond to common correlation structures, such as the
exchangeable and AR-1. 

The GEE approach yields a consistent estimator of $\mbt$ even when the
working correlation structure $\bR(\mba)$ is misspecified. However,
under such misspecification, the GEE estimator of $\mbt$ is not
efficient. Furthermore, \cite{crowder:95} established that there are
difficulties with estimating the nuisance parameter vector $\mba$ in
the GEE framework and that in certain cases, the estimator of $\mba$
does not exist.

\cite{hansen:82} proposed the class of {\em generalized method of
moments} (GMMs) for estimating the vector of regression parameters
from a set of score functions, where the dimension of the score
function exceeds that of the regression parameter.  Hansen established
that, under certain regularity conditions, the GMM estimator is
consistent, asymptotically normal and asymptotically efficient. To
overcome the difficulties associated with the GEEs, \cite{qu:00}
(henceforth abbreviated QLL) applied the principle of the GMMs in the
GEE framework that implicitly estimates the underlying correlation
structure: In particular, they proposed a clever approach based on the
{\em quadratic inference functions} (QIFs) (1) to estimate the working
correlation structure $\bR(\mba)$ such that the resulting estimator
always exists and (2) to obtain better efficiency in estimating $\mbt$
within the assumed family even under the misspecification of
$\bR(\mba)$.

\subsection{Main Results}

In this article, we derive a unified large-sample theoretical
framework for the QIF. The cornerstone of our theory is stated in
Theorem \ref{lem:quadap} which establishes a uniform quadratic
approximation to the QIF surface in a neighborhood of $\mbt_0$, the
true regression parameter vector.  This result has two fundamental
consequences. First, it provides the necessary machinery to establish
that the QIF is asymptotically equivalent to the generalized least
squares criterion. This leads to techniques for deriving large-sample
results for the QIF estimators and test statistics, analogous to
standard inferential theory for the generalized least squares methods.
Second, it provides a flexible algorithm for finding the QIF
estimators. Building on the quadratic approximation to the QIF, we
create an iteratively reweighted generalized least squares (IRGLS)
algorithm for estimation and testing in the QIF framework. This
algorithm is stable and computationally more feasible than the
Newton-Raphson algorithm recommended in the existing QIF literature
(Figure \ref{fig:qifsurf} demonstrates the necessity for the IRGLS
algorithm).

The QLL article has made important contributions to the analysis of
correlated data. However, some of the asymptotic results in QLL are
flawed, including the proof of their Theorem 1.  QLL employ a Taylor
series expansion in proving the Theorem 1; however, there is a
fundamental problem since the Taylor series expansion is in terms of
the partitioned parameter vector \(\mbt = (\mbp, \mbl)\) and not in
terms of the sample size \(N\). Even with the remainder terms,
Taylor's theorem does not provide any statements regarding the
behavior of the error term as a function of the sample size $N$ and
neither has this been established in QLL. It is not possible to patch
up either the original proof or other asymptotic results in the QLL
article (the reasons will be made clear in later sections). These
difficulties have not been identified explicitly; hence, there is not
yet any rigorous development of a unified asymptotic theory for the
QIF. We develop a more broadly applicable inferential theoretical
framework for the QIF that extends and corrects the results of QLL.

As a first step in achieving the objectives, we re-formulate the
estimated covariance matrix of the ``extended score functions''
defining the QIF. This formulation yields a consistent estimator of
the covariance matrix, which in turn lays the foundation for deriving
valid inferential theory for the QIF. The main results are summarized
as follows: 
\bee
\vspace{-0.1in}
\item We formulate an appropriate estimated covariance matrix for the
   set of extended score functions defining the QIF (Section
   \ref{sec-qif}). In Sections \ref{sec-sec2} and \ref{sec-Qresults},
   we first formulate a unified large-sample theoretical framework for
   the QIF and next derive several important asymptotic properties for
   the QIF. These lay the necessary foundation for the development of
   the asymptotic results derived in later part of the article.
\vspace{-0.1in}
\item In Section \ref{sec-test}, we first derive the principal
   result, the quadratic approximation to the QIF surface in a
   neighborhood of $\mbt_0$, the true regression parameter
   vector. Next, we formulate a statistic based on the QIF for testing
   general linear hypotheses involving correlated data. Building on
   the quadratic approximation to the QIF, we establish the asymptotic
   distribution of the generalized QIF test statistic under both the
   null and local alternative hypotheses.
\vspace{-0.1in}
\item In Section \ref{sec-IRGLS}, we propose a stable and
   computationally feasible IRGLS algorithm for estimating \(\mbt\) in
   the QIF framework. This algorithm is a step in the direction of
   developing a unified framework for estimation, testing and model
   selection for correlated data within the QIF setting.
\vspace{-0.1in}
\item In Section \ref{sec-example}, we illustrate the methods using a
   benchmark dataset consisting of the correlated binary data
   measuring the respiratory health effects of indoor and outdoor air
   pollution.  
\vspace{-0.1in}
\item In Section \ref{sec-sim}, we investigate the effect of the basis
   matrices (defining the set of extended score functions) on the size
   and power of the QIF test through Monte Carlo simulation
   experiments from Bernoulli and Gaussian distributions.
\eee
This article derives asymptotic theory for testing general linear
hypotheses based on the quadratic approximation of the QIF. However, a
common thread underlying the recent literature in the context of
nonlinear testing problems is in fact the Theorem \ref{lem:quadap}:
(1) \cite{pilla:05} derived asymptotic distribution of the test
statistic for order-restricted hypothesis testing problem; (2)
\cite{pilla:04} developed inferential theory for testing under the
general convex cone alternatives for correlated data; and (3)
\cite{loader1:05} derived several properties of the IRGLS
algorithm which is more generally applicable while providing a
flexible technique for estimation, testing and model selection with
correlated data.

Having the correct asymptotic theory for the QIF is essential for
further extensions as well as applications of the QIF, especially
given that the QIF is elegant, simple and practical to implement with
the proposed IRGLS algorithm for the analysis of correlated data.

\subsection{The Quadratic Inference Function}
\label{sec-qif}

QLL showed that the principle of GMMs can be applied in the GEE
framework by implicitly estimating the underlying correlation
structure. In particular, they assumed that the inverse of the working
correlation matrix $\bR(\mba)$ can be expressed as a linear
combination of pre-specified basis matrices \(\bM_1, \ldots, \bM_s\)
such that
\bea
  \label{eq:linearinv}
  \bR^{-1}(\mba) = \sum_{j = 1}^s \alpha_j \, \bM_j, 
\eea
where $\alpha_1, \ldots, \alpha_s$ are unknown constants. For this
article, we choose \(\bM_1\) as the identity (of appropriate
dimension) and \(\bM_2, \ldots, \bM_s\) according to the form of the
assumed underlying correlation structure. For example, (1) if
$\bR(\mba)$ is an exchangeable correlation matrix, then \(s = 2\) and
\(\bM_2\) is a matrix of 1s; and (2) if $\bR(\mba)$ is an AR-1
correlation matrix, then \(s = 3\), $\bM_2$ takes 1 on the two main
off diagonals and zero elsewhere and $\bM_3$ takes 1 at the elements
$(1, 1)$ and $(n, n)$ and zero elsewhere.  In general, $\bM_3$ is a
minor boundary correction and can be omitted.  If the covariate is
time-independent, then the boundary correction does not have an effect
on the inference since the corresponding components of the score
vector (see equation (\ref{eq:gfn} below) are linearly dependent on
other terms.

The quasi-score estimating equations in (\ref{eq:quasi}) can be
expressed, under the representation (\ref{eq:linearinv}), as 
\bes
  \sum_{i = 1}^N \nabla \bh_i^{T} \, \bA_i^{-1/2} \; (\alpha_1 \,
  \bM_1 + \cdots + \alpha_s \, \bM_s) \, \bA_i^{-1/2} \, (\bY_i - 
  \bh_i). 
\ees
These estimating equations are linear combinations of
elements of a set of {\em extended score functions}
\begin{eqnarray*}
  \overline{\bg}_N(\mbt) &:=& \frac{1}{N} \sum_{i = 1}^N \bg_i(\mbt),
\end{eqnarray*}
where a set of ``subject-specific'' basic score functions is defined
as 
\begin{eqnarray}
  \label{eq:gfn} 
  \bg_i(\mbt) = \left\{ \begin{array}{c} \nabla \bh_i^{T} \;
  \bA_i^{-1/2} \; \bM_1 \; \bA_i^{-1/2} \; (\bY_i - \bh_i) \\
  \vdots \\ \nabla \bh_i^{T} \; \bA_i^{-1/2} \; \bM_s \;
  \bA_i^{-1/2} \, (\bY_i - \bh_i) \end{array} \right\} \quad
  \mbox{for} \quad i = 1, \ldots, N.
\end{eqnarray}

Instead of directly estimating the parameters \(\alpha_1, \ldots,
\alpha_s\), the QIF introduces a sample covariance matrix in order to
combine the score functions in an optimal manner. In general, the
equation \(\overline{\bg}_N(\mbt) = 0\) has no solution, since its
dimension is greater than the number of unknown parameters. Instead,
the parameter vector $\mbt$ is estimated by minimizing the QIF defined
as
\begin{eqnarray}
  \label{eq:qif}
  Q_N(\mbt) := N \; \overline{\bg}_N^{T}(\mbt) \; \wC_N^{-1}(\mbt)
  \; \overline{\bg}_N(\mbt),   
\end{eqnarray}
where an estimator of the {\em second moment matrix} of $\bg_1(\mbt)$
is 
\begin{equation}
  \label{eq:cN}
  \wC_N(\mbt) := \frac{1}{N} \sum_{i = 1}^N \bg_i(\mbt) \;
  \bg_i^{T}(\mbt).
\end{equation}
If the extended score vector \(\overline{\bg}_N(\mbt)\) has mean zero,
then \(N^{-1} \, \wC_N(\mbt)\) is an estimator of the covariance
matrix of \(\overline{\bg}_N(\mbt)\). The function \(Q_N(\mbt)\)
measures the size of the score vector relative to its covariance
matrix and large values of \(Q_N(\mbt)\) can be considered as an
evidence against a particular value of \(\mbt\). In this sense, the
QIF plays a role similar to the negative of the loglikelihood in
parametric statistical inference.  In particular, one can construct a
goodness-of-fit test statistic $Q_N(\mbt)$ for testing the model
assumption in (\ref{eq:link}). It follows from the results of
\cite{hansen:82} that the asymptotic distribution of $Q_N(\wmbt)$
is $\chi^2$ with $\{{\rm dim}(\bg) - {\rm dim}(\mbt)\}$ degrees of
freedom under the model assumption.

The second moment matrix estimator defined in (\ref{eq:cN}) is an
average and hence under certain regularity conditions it will converge
to the true covariance of \(\bg_1(\mbt)\) as \(N \to
\infty\). This convergence result is fundamental to adapting the
large-sample framework of the GMMs \citep{hansen:82} to the QIF
setting. The role of the matrix \(\wC_N^{1/2}(\mbt)\) is similar to
that of the matrix ${\bf a}_N^*$ defined on p.\ 1040 of
\cite{hansen:82}, which is also required to converge to a
non-degenerate limit.

Our covariance estimator $\wC_N(\mbt)$ differs from that of QLL,
who define a covariance \(\bC_N\) with a
factor of \(N^{-2}\) and correspondingly omit the factor of \(N\)
from the QIF in (\ref{eq:qif}). This has led to a number of imprecise
claims in QLL, centered around their statement on p.\ 829 that
\(\bC_N\) converges to \(E(\bC_N)\). In fact, their \(\bC_N\)
converges to zero. Furthermore, the asymptotic result for $\wmbt_N$ in
Section 3.6 of QLL is incorrect since $\wC_N$ in equation (8) of their
article is approaching zero as $N \to \infty$ and is not a consistent
estimator of $\bS_{\mbt_0}(\mbt)$, the true covariance matrix of
$\bg_1(\mbt)$ evaluated at $\mbt_0$. However, our definition of the
QIF matches that presented by \cite{park:99}.

While the correlation model (\ref{eq:linearinv}) motivates our choice
of the score vector, the fundamental property \(E\{\bg_1(\mbt_0)\} =
0\) holds whether or not the covariance assumption is
correct. Similarly, \(\widehat \bC_N(\mbt_0)\) consistently estimates
\(N \, \cov\{ \overline{\bg}_N(\mbt_0) \}\). Therefore, inference
based on the QIF is semiparametric, in the sense that procedures are
asymptotically valid whether or not the covariance model is
correct. \cite{bickel:98} and \cite{kosorok:06} present a detailed
exposition of the mathematical aspects of semiparametric inference. 

\section{Large-Sample Properties of the Extended Score Functions}
\label{sec-sec2}

In order to establish the asymptotic results in any regression
problem, one must first state assumptions regarding the behavior of
the design matrices as the sample size \(N\) increases. Without
formulating such assumptions, it is not possible to establish even the
consistency of the QIF estimators. However, such assumptions are
missing from the earlier QIF work.

The main requirement for our asymptotic theory is to be able to apply
the strong law of large numbers to show that
\(\overline{\bg}_N(\mbt)\), \(\widehat{\bC}_N(\mbt)\) and other 
averages converge to appropriate non-degenerate limits. A sufficient
condition is the following ``random design'' assumption.

\begin{assumption}
\label{as1}
The pairs $(\bY_i, \bX_i^T)$ are assumed to be an independent sample
from a $\{n \times (q + 1)\}$-dimensional distribution $\bF$, where
\(\bX_i = (\bX_{i1}, \ldots, \bX_{in})\) is the $(q \times
n)$-dimensional design matrix for the \(i\)th ($i = 1, \ldots, N$)
subject.  
\end{assumption}

\begin{remark} The independence part of Assumption \ref{as1} is
between different subjects, or with respect to the index \(i\). The
elements of \(\bX_i\) need not be independent of each other; hence,
this assumption incorporates both time-dependent and time-independent
covariates. Note that there exists a dependence of \(\bY_i\) on
\(\bX_i\) through the link and variance functions in (\ref{eq:link})
and (\ref{eq:var}), respectively. 
\end{remark} 

\vspace{0.1in}
All throughout this article, $E_{\mbt_0}(\cdot)$ denotes the
expectation operator with respect to the true regression parameter
vector $\mbt_0$. Our results are based on the implicit assumption that
all expectations are finite and the convergence statements are uniform
for \(\mbt\) in bounded sets. Uniformity results for the strong law of
large numbers are derived by \cite{rubin:56}.

\begin{theorem}[Asymptotic normality of $\overline{\bg}_N(\mbt_0)$]
  \label{thm1} Let $\mbt_0$ be the true parameter vector. Under
  Assumption \ref{as1}, 
\begin{eqnarray*}
  \overline{\bg}_N(\mbt) &\overset{a.s.}{\longrightarrow}& E_{\mbt_0}
  \{ \bg_1(\mbt)\} = 0 \quad \mbox{if} \quad \mbt = \mbt_0
\end{eqnarray*}
and
\begin{equation}
  \label{eq:gn-clt}
  \sqrt{N} \, \overline{\bg}_N(\mbt_0) \overset{d}{\longrightarrow}
  N_r\{0, \bS_{\mbt_0}(\mbt_0)\},
\end{equation}
where $r = q \, s$ and \(\bS_{\mbt_0}(\mbt_0)\) is the true
covariance matrix of \(\bg_1(\mbt)\) evaluated at $\mbt_0$.
\end{theorem}

\vspace{0.2in}
It is easy to verify that \(E_{\mbt_0}\{\bg_1(\mbt_0)\} = 0\). The
following identifiability assumption is required to develop the
large-sample theory for the QIF.

\begin{assumption}
\label{as2}
The parameter \(\mbt\) is estimable, in the sense that
$E_{\mbt_0}\{\bg_1(\mbt)\} = 0$ if and only if $\mbt \neq \mbt_0$.
\end{assumption}

Another application of the strong law of large numbers establishes
that \(\wC_N(\mbt)\) converges to its expected value, a
non-degenerate limit, which is required to invoke the results of
\cite{hansen:82}. 

\begin{theorem}[Consistency of $\wC_N(\mbt)$]
  \label{thm-cn}
Under Assumptions \ref{as1} and \ref{as2},
\begin{equation}
  \label{eq:sigma}
  \wC_N(\mbt)  \overset{a.s.}{\longrightarrow} E_{\mbt_0} \left\{ 
  \bg_1(\mbt) \; \bg_1^{T}(\mbt) \right\}
  := \bS_{\mbt_0}(\mbt) \quad \mbox{as} \quad N \to \infty.
\end{equation}
\end{theorem}

In this article we follow the prescription of Assumption 3.6 of
\cite{hansen:82} which is restated for the current framework.

\begin{assumption}
\label{as3}
The matrix $\bS_{\mbt_0}(\mbt)$ is strictly positive definite.
\end{assumption}

\begin{remark} The estimator of second moment matrix
\(\wC_N(\mbt)\) may be singular. However, any vector in the null
space of \(\wC_N(\mbt)\) must be orthogonal to each of the
subject-specific score functions $\bg_i(\mbt) \, (i = 1, \ldots, N)$
and consequently to \(\overline{\bg}_N(\mbt)\). As a result, one can
replace \(\wC_N^{-1}(\mbt)\) by any generalized inverse such as the
Moore-Penrose generalized inverse.
\end{remark} 

\vspace{0.1in}
We first state several important assumptions required to establish the
large-sample properties of $\overline{\bg}_N(\mbt)$.

\begin{assumption}
\label{as4} The parameter space of $\mbt$ denoted by $\mB$ is 
compact.
\end{assumption}

\begin{assumption}
\label{as6}
The expectation $E_{\mbt_0}\{\overline{\bg}_N(\mbt)\}$ exists and is
finite for all $\mbt \in \mB$ and is continuous in $\mbt$.
\end{assumption}

The compactness assumption is necessary to invoke the uniformity
results of \cite{rubin:56} and it is unavoidable since the QIF surface
is often not convex. For non-compact parameter spaces, Lemma
\ref{lemma:bn-con} is only applicable to a sequence of local minima.

From Theorem 2.1 of \cite{hansen:82}, the following result holds. 

\begin{lemma}[Consistency of $\wmbt_N$]
  \label{lemma:bn-con}
  Under Assumptions \ref{as4}--\ref{as6}, the QIF estimator 
\bes
  \wmbt_N := \underset{\mbt \, \in \, \mB}{\arg \, \min} \; Q_N(\mbt)
\ees
exists and $\wmbt_N \overset{a.s.}{\longrightarrow} \mbt_0$ as \(N \to
  \infty\).  
\end{lemma}

\vspace{0.2in} 
Our next goal is to derive the asymptotic distribution of
$\wmbt_N$. Let 
\begin{eqnarray*}
  \label{eq:d0}
  \bD(\mbt) := E_{\mbt_0} \left\{ \frac{\partial}{\partial
  \mbt} \, \bg_1(\mbt) \right\} = E_{\mbt_0} \{ \nabla 
  \bg_1(\mbt)\}.  
\end{eqnarray*}
Once again, from the strong law of large numbers, it follows that 
\bea
  \label{eq:dotg}
  \nabla \overline{\bg}_N(\mbt) \overset{a.s.}{\longrightarrow}
  E_{\mbt_0}\{ \nabla \bg_1(\mbt)\} = \bD(\mbt).
\eea

The extended score vector $\overline{\bg}_N(\mbt)$ is a random vector
and hence $\nabla \, \overline{\bg}_N(\mbt)$ is a random
matrix. Therefore, the claims on p.\ 829 of QLL that $\nabla \,
\overline{\bg}_N(\mbt)$ is nonrandom and $E\{\nabla
\overline{\bg}_N(\mbt) \} = \overline{\bg}_N(\mbt)$ cannot be true. 

\begin{assumption}
\label{as7}
The subject-specific score functions \(\bg_i(\mbt)\) $(i = 1, \ldots,
N)$ have uniformly continuous second-order partial derivatives with
respect to the elements of the vector \(\mbt\).
\end{assumption}

Owing to Theorems 3.1 and 3.2 of \cite{hansen:82}, the following
result holds.

\begin{theorem}[Asymptotic normality of $\wmbt_N$]
Under Assumptions \ref{as1}--\ref{as7}, the asymptotic distribution
of $\wmbt_N$ is 
\bea
  \label{eq:basymp}
  \sqrt{N} \left(\wmbt_N - \mbt_0 \right) \overset{d}{\longrightarrow} 
  N_q \left\{0, \bJ^{-1}(\mbt_0)\right\}
\eea
as \(N \to \infty\), where
\bea
  \label{eq:J}
  \bJ(\mbt_0) = \bD^T(\mbt_0) \; \bS^{-1}_{\mbt_0}(\mbt_0) \;
  \bD(\mbt_0). 
\eea
\end{theorem}

The claim by QLL on p.\ 835, below equation (A1), that $(\wmbt_N -
\mbt_0)$ converges in law to a normal distribution is incorrect. In
fact, \(\sqrt{N}(\wmbt_N - \mbt_0)\) has an asymptotic normal
distribution while \((\wmbt_N - \mbt_0)\) converges in probability to
zero. The matrix \(\bd_0\) is not defined and it is not possible to
define this in a manner consistent with the remainder of their
article. While the statement of Theorem 1 in QLL requires the
non-centrality parameter (consequently the partitioned matrices
\(\bJ_{\mbp \mbp}\) etc., \(\bd_0\) and \(\bS\)) to be \(O(1)\), the 
equation (A1) on p.\ 835 requires \(\bJ_{\mbp \mbp}\) to be
\(O_p(N^{-1})\). We believe, based on comparison with other work on
the QIF, that the authors probably intended to write \(\bd_0 =
E\{\nabla \overline{\bg}_N(\mbt)\} = O_p(1)\) and \(\bS =
\hbox{cov}\{\overline{\bg}_N(\mbt) \} = O(N^{-1})\). However, this
means that the statement of their Theorem 1 is incorrect.  

\section{Fundamental Results for the Quadratic Inference Functions}
\label{sec-Qresults}

In this section, we establish several fundamental results for the QIF
which lay the foundation for deriving the asymptotic distribution of
inference functions presented in the next section.

One main focus is on the vector \(\nabla Q_N(\mbt)\) of partial
derivatives and matrix \(\nabla^2 Q_N(\mbt)\) of second-order partial
derivatives.  Along the way, we derive the correct versions of several
claims made by QLL.  For example, on p.\ 830, QLL incorrectly claim
that \(\nabla^2 Q_N(\mbt)\) converges in probability. In fact, the
second derivative matrix has size \(O_p(N)\).

\begin{assumption}
\label{as8}
The first and second-order partial derivatives of
\(\overline{\bg}_N(\mbt)\) and \(\wC_N(\mbt)\) have finite means. 
\end{assumption}

\begin{theorem}
  \label{thm:qn-clt}
Under Assumption \ref{as8}, 
\begin{equation}
  \label{eq:q-clt}
  \frac{1}{2 \, \sqrt{N}} \; \nabla Q_N(\mbt_0)
  \overset{d}{\longrightarrow} N_q \left\{0, \bJ(\mbt_0)\right\}
\end{equation}
as \(N \to \infty\). There exists a non-random matrix $V(\mbt)$,
continuous in $\mbt$, such that
\begin{equation}
  \label{eq:qd-clt}
  \frac{1}{2N} \nabla^2 Q_N(\mbt) = \bV(\mbt) + o_p(1),
\end{equation}
where $\bV(\mbt_0) = \bJ(\mbt_0)$ and \(o_p(1)\) error term is uniform
on compact sets.  
\end{theorem}

{\em Proof:} Differentiating $Q_N(\mbt)$, with respect to the
\(k\)th ($k = 1, \ldots, q$) element \(\beta_k\) of \(\mbt\) yields
\begin{eqnarray}
  \frac{\partial}{\partial \beta_k} \, Q_N(\mbt)
     &=& 2 \, N \; \overline{\bg}_N^{T}(\mbt) \; \wC_N^{-1}(\mbt) \;
       \frac{\partial \, \overline{\bg}_N(\mbt)}{\partial \beta_k}
     + N \; \overline{\bg}_N^{T}(\mbt) \;
       \frac{\partial \, \wC_N^{-1}(\mbt)}{\partial \beta_k}
       \; \overline{\bg}_N(\mbt). \nonumber \\
  \label{eq:difqn}
\end{eqnarray}
From Theorem \ref{thm1}, it follows that at \(\mbt = \mbt_0\),
\(\sqrt{N} \; \overline{\bg}_N(\mbt_0) = O_p(1)\) and \(\partial \,
\wC_N^{-1}(\mbt)/\partial \beta_k\) has a finite limit by the strong
law of large numbers. Consequently, the second term in
(\ref{eq:difqn}) is \(O_p(1)\). Therefore, it follows that
\begin{eqnarray*}
  \frac{1}{2 \, \sqrt{N}} \, \nabla {Q}_N(\mbt_0) &=& \sqrt{N} \; 
  \nabla \overline{\bg}_N^{T}(\mbt_0) \; \wC_N^{-1}(\mbt_0) \;
  \overline{\bg}_N(\mbt_0) + o_p(1) \\
  &=& \sqrt{N} \; \bD^T(\mbt_0) \; \bS^{-1}_{\mbt_0}(\mbt_0) \;
  \overline{\bg}_N(\mbt_0) + o_p(1).
\end{eqnarray*}
The last equation follows from the result (\ref{eq:sigma}) in Theorem
\ref{thm1} and the result (\ref{eq:dotg}). The asymptotic distribution
of \(\sqrt{N} \, \overline{\bg}_N(\mbt_0)\) given in (\ref{eq:gn-clt})
yields the result (\ref{eq:q-clt}).

Similarly,
\begin{eqnarray*}
  \frac{1}{N} \, \frac{\partial^2}{\partial \beta_j \, \partial
  \beta_k}  \, Q_N(\mbt) 
  &=& 2 \; \overline{\bg}_N^{T}(\mbt) \; \wC_N^{-1}(\mbt) \;
  \frac{\partial^2 \, \overline{\bg}_N(\mbt)}{\partial \beta_j \,
  \partial \beta_k} + \overline{\bg}_N^{T}(\mbt) \; \frac{\partial^2
  \, \wC_N^{-1}(\mbt)}{\partial \beta_j \, \partial \beta_k} \;
  \overline{\bg}_N(\mbt) \\ 
   && + \; 2 \; \frac{\partial \, \overline{\bg}_N^{T}(\mbt)}{\partial
   \beta_j} \; \wC_N^{-1}(\mbt) \; \frac{\partial \,
   \overline{\bg}_N(\mbt)}{\partial \beta_k} + 2  \; \frac{\partial \,
   \overline{\bg}_N^{T}(\mbt)}{\partial \beta_j} \; \frac{\partial \,
   \wC_N^{-1}(\mbt) }{\partial \beta_k} \; \overline{\bg}_N(\mbt) \\ 
   && + \; 2  \; \frac{\partial \,
   \overline{\bg}_N^{T}(\mbt)}{\partial \beta_k} \; 
         \frac{\partial \, \wC_N^{-1}(\mbt)}{\partial \beta_j} \;
	 \overline{\bg}_N(\mbt). 
\end{eqnarray*}
As $N \to \infty$, by the strong law of large numbers, each of these
terms is converging to a non-degenerate limit, which in turn lead to
result (\ref{eq:qd-clt}) for an appropriate \(\bV(\mbt)\).  At \(\mbt
= \mbt_0\), all the terms involving \(\overline{\bg}_N(\mbt)\)
converge to zero, leaving only the third term:
\[
  \frac{1}{N} \frac{\partial^2 Q_N(\mbt_0)}{\partial \beta_j \,
  \partial \beta_k} = 2 \; \frac{\partial \, 
  \overline{\bg}_N^{T}(\mbt_0)}{ \partial \beta_j} \;
  \wC_N^{-1}(\mbt_0) \; \frac{ \partial \, 
  \overline{\bg}_N(\mbt_0)}{\partial \beta_k} + o_p(1).
\]
In matrix form, this can be restated as 
\begin{equation}
  \label{eq:qddot}
  \frac{1}{2N} \, \nabla^2 Q_N(\mbt_0) = \bJ(\mbt_0) + o_p(1)
\end{equation}
which establishes the result (\ref{eq:qd-clt}) at $\mbt = \mbt_0$,
implying that \(\bV(\mbt_0) = \bJ(\mbt_0)\). 
\hfill \rule{2mm}{2mm}

From Theorem \ref{thm:qn-clt}, we have the following corollary.

\begin{corollary} 
\label{cor:Qderiv} The first and second derivatives of  $Q_N(\mbt)$
at $\mbt = \mbt_0$ satisfy $\nabla Q_N(\mbt_0) = O_p(\sqrt{N})$ and
$\nabla^2 Q_N(\mbt_0) = O_p(N)$, respectively.
\end{corollary}

The analysis of the first derivative vector as well as the second
derivative matrix of $Q_N(\mbt)$ are critical to the development of
asymptotic theory for the QIF and numerical algorithms to find the
QIF estimators. However, several inadequacies exist in the previous
results provided by QLL. First, their derivation involves
multiplication of three and four-way arrays which are not clearly
defined, leading to an incorrect expression for the second derivative
matrix, $\nabla^2 Q_N(\mbt_0)$. The claims made by QLL about
convergence of the first and second derivatives of $Q_N(\mbt)$
contradict the results shown in Corollary \ref{cor:Qderiv}.  These
problems in QLL lead to their claim on p.\ 835 that \((\widehat{\mbp}
- \mbp_0)\) and \((\widehat{\mbl} - \mbl_0)\) have a limiting normal
distribution, however, with a missing a factor of \(\sqrt{N}\).

\section{Testing for General Linear Hypotheses within the QIF
Framework} 
\label{sec-test}

In this section we first establish the quadratic approximation of the
QIF in a local neighborhood of \(\mbt_0\). Next, we derive an
asymptotic distribution of the test based on the QIF for testing a
general linear hypothesis and demonstrate that the Theorem 1 of QLL
becomes a special case of our result.

\subsection{Asymptotic Distribution of the Inference Functions}

The fundamental principle underlying the large-sample results
presented in this article is a quadratic approximation of the QIF in a
local neighborhood of the true regression parameter vector
\(\mbt_0\). Under this approximation, the problem of minimizing the
QIF is asymptotically equivalent to a generalized least squares
criterion. Consequently, standard asymptotic results from linear
models can be applied to the QIF framework.

\noindent{\em Definition:} A ball of radius $o(\sqrt{N})$ is defined
as $\{\mbx: \| \mbx\| \leq r_N \}$, where $\{r_N\!: N > 1\}$ is a
sequence of constants with \(r_N = o(\sqrt{N})\).

For exposition, we define 
\bea
  \label{eq:ZN}
  \bZ_N := \frac{1}{2 \, \sqrt{N}} \; \nabla Q_N(\mbt_0).
\eea

\begin{theorem}[Quadratic approximation of the QIF]
\label{lem:quadap}
For a fixed $q$-dimensional vector \(\mbx\), the following
representation holds:
\begin{equation}
  Q_N \left(\mbt_0 + N^{-\frac{1}{2}} \, \mbx \right) = Q_N(\mbt_0)
    + 2 \; \ip{\mbx}{\bZ_N} + \mbx^{T} \; \bJ(\mbt_0) \; \mbx + o_p(1)
    \label{eq:qnlocal} 
\end{equation}
as \(N \to \infty\), where $\ip{\cdot}{\cdot}$ is the vector inner
product and the \(o_p(1)\) term is uniform for \(\mbx\) in a ball of
radius \(r_N = o(\sqrt{N})\).   
\end{theorem}

\emph{Proof:} The Taylor series expansion yields
\bes
  Q_N \left(\mbt_0 + N^{-\frac{1}{2}} \, \mbx \right) = Q_N(\mbt_0)
    + 2 \; \ip{\mbx}{\bZ_N} + \frac{1}{2N} \; \mbx^{T} \; \nabla^2
    Q_N \left(\mbt_N^{\dag} \right) \; \mbx,
\ees
where \(\mbt_N^{\dag}\) lies between \(\mbt_0\) and \( (\mbt_0 +
N^{-1/2} \, \mbx)\) for each \(N\). From the uniformity result
(\ref{eq:qd-clt}), it follows that   
\bes
  \frac{1}{2N} \, \nabla^2 Q_N \left(\mbt_N^{\dag} \right) =
    \bV\left(\mbt_N^{\dag} \right) + o_p(1).
\ees
As \(N \to \infty\), 
\(\mbt_N^{\dag} \longrightarrow \mbt_0\). From the continuity result
of Theorem \ref{thm:qn-clt}, it follows that \(\bV(\mbt_N^{\dag})
\longrightarrow \bV(\mbt_0) = \bJ(\mbt_0)\), yielding the desired
result. \hfill \rule{2mm}{2mm} 

\begin{corollary} The quadratic approximation in Theorem
\ref{lem:quadap} can be expressed as
\begin{eqnarray*}
    Q_N \left(\mbt_0 + N^{-\frac{1}{2}} \, \mbx \right) &=&
    \left\{ \bZ_N + \bJ(\mbt_0) \, \mbx \right\}^T \; \bJ^{-1}(\mbt_0)
    \; \left\{ \bZ_N + \bJ(\mbt_0) \, \mbx \right\} \\
   && + \; Q_N(\mbt_0) - \bZ_N^T \; \bJ^{-1}(\mbt_0) \; \bZ_N + 
    o_p(1). 
\end{eqnarray*}
\end{corollary}

The representation in the above corollary establishes that the QIF is
asymptotically equivalent to a generalized least squares
criterion. This simplifies derivation of large-sample results, since
known properties of the weighted least squares will hold
asymptotically for the QIF. Second, it leads to an IRGLS algorithm for
finding the QIF estimator \(\wmbt_N\).

The minimizer \(\mbx_N^{\star}\) of the quadratic approximation in
(\ref{eq:qnlocal}) is given by 
\bes
  \mbx_N^{\star} = -\bJ^{-1}(\mbt_0) \, \bZ_N.
\ees
Since \(\bZ_N\) has a limiting distribution, it follows that
\(\mbx_N^{\star}\) lies in the ball of radius \(r_N\) with probability
converging to 1. This result, combined with the uniformity of the
error term in (\ref{eq:qnlocal}), yields the following corollaries.

\begin{corollary} Let \(\widehat{\mbx}_N\) be the minimizer of
\(Q_N\left(\mbt_0 + N^{-1/2} \, \mbx \right)\), then \(\wmbt_N =
(\mbt_0 + N^{-1/2} \, \widehat{\mbx}_N)\). Equivalently, 
\bes
  \widehat{\mbx}_N = -\bJ^{-1}(\mbt_0) \, \bZ_N + o_p(1)
\ees
and
\begin{equation}
  \label{eq:bhatp}
  \wmbt_N = \mbt_0 - N^{-\frac{1}{2}} \, \bJ^{-1}(\mbt_0) \,
  \bZ_N + o_p(N^{-1/2}).
\end{equation}
\end{corollary}
\begin{corollary} The asymptotic distribution of \(\wmbt_N\) is
given by
\bes
  \sqrt{N} \left(\wmbt_N - \mbt_0 \right) \overset{d}{\longrightarrow}
  N_q \left\{0, \bJ^{-1}(\mbt_0) \right\} \quad \mbox{as} N \to
  \infty. 
\ees
\end{corollary}
\begin{corollary} The following result holds:
\begin{equation} 
  Q_N(\wmbt_N) = Q_N(\mbt_0) - \bZ_N^{T} \; \bJ^{-1}(\mbt_0) \;
  \bZ_N + o_p(1) \quad \mbox{as} \quad N \to \infty.
  \label{eq:qnalt} 
\end{equation}
\end{corollary}
\vspace{-0.2in}

\subsection{Asymptotic Distribution of a Generalized QIF Test
Statistic} 

In this section, we derive the asymptotic theory for testing a general
linear hypotheses. Consequently, the one presented in QLL becomes a
special case. 

Following the notation in \cite{chris:02}, we consider the problem
of testing the general linear hypothesis
\bea
  \label{eq:hyp}
  H_0\!: \mbL^{T} \, \mbt = \bb \quad \mbox{versus} \quad 
  H_1\!: \mbL^{T} \, \mbt \ne \bb,
\eea
where, for some \(p < q\), the \((q \times p)\) matrix \(\mbL\)
imposes \(p\) linearly independent constraints on the parameter vector
\(\mbt\) and constant vector \(\bb \in \Re^p\).

The QIF-based test statistic for testing the general linear hypothesis
problem (\ref{eq:hyp}) is
\bea
   \label{eq:TN}
   T_N := Q_N(\wtb_N) - Q_N(\wmbt_N),
\eea
where the unrestricted and constrained minimizers of \(Q_N(\mbt)\)
respectively, are
\bes
  \wmbt_N &=& \underset{\mbt \, \in \, \mB}{\arg \, \min} \; Q_N(\mbt)
\ees
and
\bes
  \wtb_N &=& \underset{\mbt \, \in H_0}{\arg \, \min} \; Q_N(\mbt).
\ees

\begin{theorem}[Null asymptotic distribution of $T_N$]
\label{thm-test}
For testing the hypothesis problem (\ref{eq:hyp}), the QIF-based test
statistic has the asymptotic representation
\bea
  \label{eq:test}
  T_N = \{\mbL^{T} \; \bJ^{-1}(\mbt_0) \; \bZ_N\}^{T} \; 
        \{\mbL^{T} \; \bJ^{-1}(\mbt_0) \; \mbL\}^{-1} \; 
        \{\mbL^{T} \; \bJ^{-1}(\mbt_0) \; \bZ_N\} + o_p(1),
\eea
where $\bJ(\mbt_0)$ is defined in (\ref{eq:J}).  The asymptotic
distribution of \(\, T_N\) under the null hypothesis $H_0$ is
\bes
  T_N \overset{d}{\rightarrow} \chi^2_p \quad \mbox{as} \quad N \to
  \infty, 
\ees
where $p$ is the number of linearly dependent constraints imposed by
the matrix \(\mbL\). 
\end{theorem}

\emph{Proof:} Suppose that the null hypothesis $H_0$ is true, then 
\(\mbL^{T} \, \mbt_0 = \bb\). Note that $\mbx = \sqrt{N} \, (\mbt -
\mbt_0)$ implies $\mbL^T \, \mbx = \sqrt{N} (\mbL^T \, \mbt - 
\mbL^T \, \mbt_0) = \sqrt{N} (\mbL^T \, \mbt - \bb)$. Therefore,
minimizing the QIF in (\ref{eq:qnlocal}) subject to $\mbL^T \,
\mbt = \bb$ is equivalent to minimizing it subject to 
$\mbL^T \, \mbx = 0$. Following arguments similar to those in the
previous section, we obtain
\bes
  \wtb_N = \mbt_0 - N^{-\frac{1}{2}} \,
  \bJ^{-1}(\mbt_0) \left[ \bI - \mbL \; \left\{
  \mbL^{T} \; \bJ^{-1}(\mbt_0) \; \mbL \right\}^{-1} \;
  \mbL^{T} \; \bJ^{-1}(\mbt_0) \right] \bZ_N + o_p(N^{-1/2}) 
\ees
and 
\begin{eqnarray}
  Q_N(\wtb_N) - Q_N(\mbt_0) &=&
      + \; \{\mbL^{T} \; \bJ^{-1}(\mbt_0) \, \bZ_N\}^{T} \; 
      \{\mbL^{T} \; \bJ^{-1}(\mbt_0) \; \mbL\}^{-1}
      \; \{\mbL^{T} \; \bJ^{-1}(\mbt_0) \, \bZ_N\} \nonumber \\
      && \; - \bZ_N^{T} \; \bJ^{-1}(\mbt_0) \; \bZ_N + o_p(1).
  \label{eq:qnnull}
\end{eqnarray}
Combining the results (\ref{eq:qnalt}) and (\ref{eq:qnnull}), we
obtain (\ref{eq:test}).  Equation (\ref{eq:q-clt}) implies
that \(\mbL^{T} \bZ_N\) has an asymptotic \(N_p\left\{0, \mbL^{T} \;
\bJ^{-1}(\mbt_0) \; \mbL\right\}\) distribution. The asymptotic
\(\chi^2_p\) distribution of \(T_N\) as \(N \to \infty\) follows
immediately. \hfill \rule{2mm}{2mm}  


\subsection{Testing Under Local Alternatives}

We consider the hypothesis testing problem (\ref{eq:hyp}), but
assume that the alternative hypothesis is true. Specifically, consider
a sequence of local alternative parameter vectors \(\mbt_N = (\mbt_0 +
N^{-1/2} \mbv)\), where \(\mbL^T \mbt_0 = \bb\) and $\mbv$ is a fixed
$q$-dimensional vector.

In order to establish the large-sample properties of the test
statistic \(T_N\) under this model, we can proceed essentially as
before with the exception that \(\sqrt{N} \,
\overline{\bg}_N(\mbt_0)\) has non-zero mean. The multivariate
central limit theorems become, respectively
\begin{eqnarray*}
  \sqrt{N} \; \overline{\bg}_N(\mbt_0) &\overset{d}{\longrightarrow}& 
  N_r\{ - \bD(\mbt_0) \, \mbv, \, \bS_{\mbt_0}(\mbt_0)\} 
\end{eqnarray*}
and
\begin{eqnarray*}
  \bZ_N = \frac{1}{2 \sqrt{N}} \; \nabla Q_N(\mbt_0)
    &\overset{d}{\longrightarrow}& N_q \left\{-\bJ(\mbt_0) \;
    \mbv, \, \bJ(\mbt_0) \right\}
\end{eqnarray*}
as \(N \to \infty\). The asymptotic representation (\ref{eq:test})
continues to hold under the local alternatives; therefore, we have the
following asymptotic distribution for \(T_N\).

\begin{theorem}[Asymptotic distribution of $T_N$ under local
alternatives] 
\label{thm:alttest}
Under $\mbt_N = (\mbt_0 + N^{-1/2} \mbv)$, the asymptotic distribution
of the test statistic $T_N$ is non-central chi-squared with a
non-centrality parameter defined as
\bea
  \label{eq:nparm}
  \delta^2 := \mbv^{T} \; \mbL \; \{\mbL^{T} \;
  \bJ^{-1}(\mbt_0) \; \mbL\}^{-1} \; \mbL^{T} \, \mbv.
\eea
\end{theorem}

\vspace{0.1in}
{\em Example:} QLL partitioned the regression parameter vector as
\(\mbt^T = (\mbp^T, \mbl^T)\) and considered testing the hypothesis
\(H_0\!: \mbp = \mbp_0\). This corresponds to $\mbL^{T} = (\bI \quad
{\bf 0})$ and \(\bb = \mbp_0\) in the hypothesis problem
(\ref{eq:hyp}). 

The result of Theorem \ref{thm:alttest} is applicable if we partition
the asymptotic covariance matrix of $\wmbt_N$ as
\bes
  \bJ(\mbt_0) = \begin{pmatrix}
  \bJ_{\mbp_0 \mbp_0} & \bJ_{\mbp_0 \mbl_0} \\
  \bJ_{\mbl_0 \mbp_0} & \bJ_{\mbl_0 \mbl_0}
  \end{pmatrix},
\ees
where $\mbt_0 = (\mbp_0, \mbl_0)$ is the null value of $\mbt =
(\mbp, \mbl)$. From standard results for the inversion of
a partitioned matrix, the non-centrality parameter can be expressed 
as
\bes
  \delta^2 = \mbv^{T} \, \mbL \; \left(\bJ_{\mbp_0 \mbp_0} -
  \bJ_{\mbp_0 \, \mbl_0} \; \bJ_{\mbl_0 \, \mbl_0}^{-1} \; \bJ_{\mbl_0
  \mbp_0} \right)^{-1} \; \mbL^{T} \, \mbv.
\ees
This agrees with the result in QLL, subject to the
concerns over the scaling of \(\bJ_{\mbp \mbp}\) discussed earlier. 

\section{The IRGLS Algorithm}
\label{sec-IRGLS}

In this section, we derive the iteratively reweighted generalized
least squares (IRGLS) algorithm for finding the QIF estimator of
\(\mbt\). The necessity for such an algorithm is illustrated first
using a simulated experimental data.

\begin{figure}
  \centerline{\scalebox{0.8}{\includegraphics{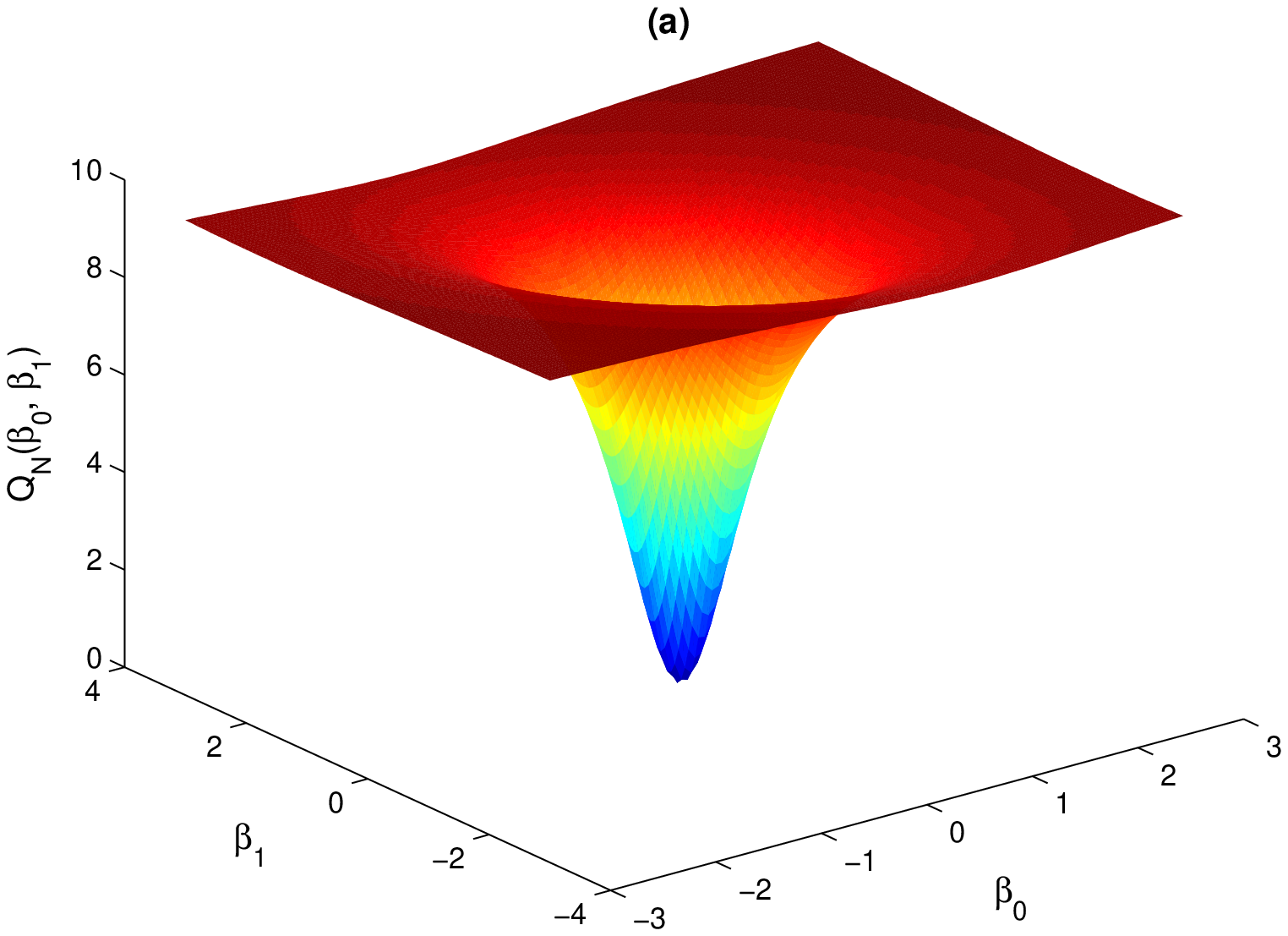}} }
  \vspace{-0.2in}
  \centerline{\scalebox{0.8}{\includegraphics{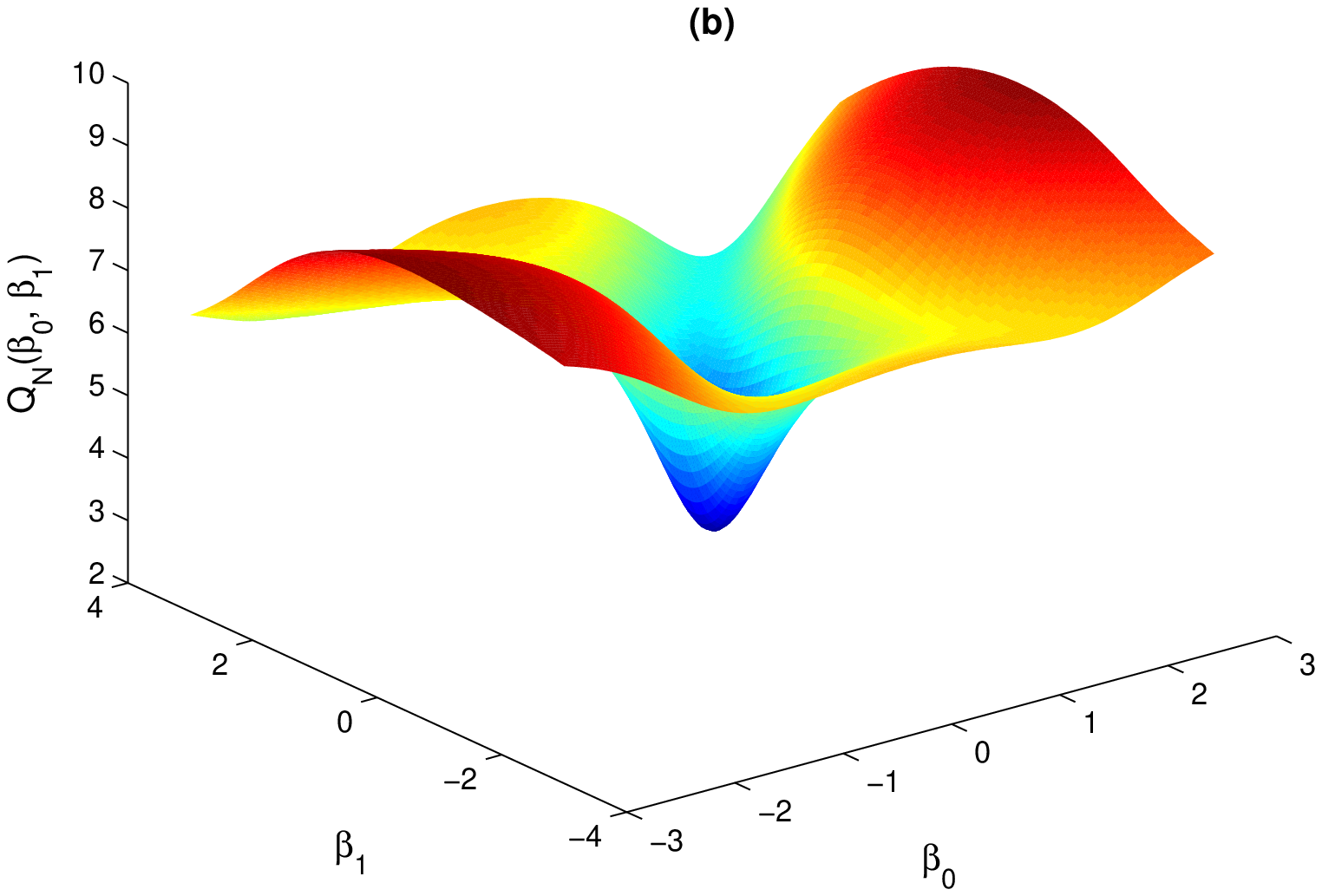}} }
  \vspace{-0.2in}
\caption{The surface plot of \(Q_N(\mbt)\), where $\mbt =
(\beta_0, \beta_1)^T$, under the AR-1 correlation structure for the
correlated (a) Gaussian responses and (b) Bernoulli responses.}
\label{fig:qifsurf}
\end{figure}

We assumed ten subjects (\(i = 1, \ldots, 10\)) and four observations
per subject (\(t = 1, \ldots, 4\)) under an AR-1 correlation structure
with autocorrelation \(\rho = 0.5\). We constructed the extended score
vector $\overline{\bg}_N(\mbt)$ using $\bM_1 = \bI$ and
\bea
  \label{eq:M2}
  \bM_2 = \begin{pmatrix}
  0 & 1 & 0 & 0 \\
  1 & 0 & 1 & 0 \\
  0 & 1 & 0 & 1 \\
  0 & 0 & 1 & 0 \\
  \end{pmatrix}
\eea
as the basis matrices. The fitted models are $\mu_{it} = \beta_0 +
\beta_1 \, (t - 2.5)$ and
\bes
  \log\left\{ \frac{\mu_{it}}{(1 - \mu_{it})} \right\} = \beta_0
  + \beta_1 \, (t - 2.5) \quad \mbox{for} \quad i = 1, \ldots, 10; t =
  1, \ldots, 4,
\ees
respectively, for the Gaussian and Bernoulli responses. Let $\mbt =
(\beta_0, \beta_1)^T$. Figure \ref{fig:qifsurf} displays the QIF
surface \(Q_N(\mbt)\) for simulated correlated data generated from the
Gaussian and Bernoulli distributions, respectively.  Notice the
strikingly different behavior under these two models. For the
correlated responses from the Gaussian distribution, the QIF surface
is bounded above by \(N = 10\) and converges to \(10\) as
\(\|\mbt\| \to \infty\) in any direction. Even in such a scenario, the
Newton-Raphson algorithm can diverge. In the Bernoulli case, the
surface has multiple ridges, valleys as well as local minima as
\(\|\beta\| \to \infty\) in some directions. It is clear that, carefully
designed algorithms are necessary to reliably find the global minimum
of $Q_N(\mbt)$.

The QIF surface plots amplify the necessity for the development of a
stable algorithm for finding the QIF estimator of $\mbt$.  The
Newton-Raphson algorithm, recommended by QLL, requires accurate
starting values to converge, especially in situations resembling that
of Figure \ref{fig:qifsurf}(b). Furthermore, in order to implement the
Newton-Raphson algorithm, we need to find the matrix
\(\nabla^2 Q_N(\mbt)\) which can be a computationally daunting task
even for small $N$.

\noindent {\em IRGLS Algorithm:} The equation (\ref{eq:bhatp}) forms
the basis of our algorithm. \\
{\em Step 1.} Start with an initial value of the parameter vector
\(\mbt^{(1)}\). \\
{\em Step 2.} Find the updated value for $\mbt$ via
\bes
  \mbt^{(j + 1)} = \mbt^{(j)} - \frac{1}{2 \, N} \;
  \widehat{\bJ}^{-1}_N \left(\mbt^{(j)} \right) \; \nabla
  Q_N \left(\mbt^{(j)} \right) \quad 
  \mbox{for $j = 1, 2, \ldots$},
\ees
where
\begin{equation}
  \widehat{\bJ}_N(\mbt) := \nabla \overline{\bg}_N^T(\mbt) \; 
     \wC_N^{-1}(\mbt) \; \nabla \overline{\bg}_N(\mbt). 
  \label{eq:jhat}
\end{equation}

If the above iterative scheme converges to a limit \(\mbt^{\infty}\),
then the limit must be a stationary point satisfying $\nabla
Q_N(\mbt^{\infty}) = 0$. An S-Plus library implementing this IRGLS 
algorithm is developed by \cite{loader2:05}.

The IRGLS algorithm proposed here inherits the standard advantages of
the IRLS methods \citep{green:84} over the Newton-Raphson algorithm:
(1) it avoids the complexity of computing \(\nabla^2 Q_N(\mbt)\); and
(2) the algorithm is guaranteed to move in a descent direction of the
QIF surface. This second property ensures that the algorithm cannot
converge to a local maximum.  With simple bounds on the step size, the
IRGLS algorithm converges to the QIF estimator from almost any
starting point \citep{loader1:05}.

\section{Analysis of Respiratory Health Effects Data}
\label{sec-example}

We analyze part of the longitudinal binary data on respiratory health
effects of indoor and outdoor air pollution in six U.S. cities
measured on 537 children at ages 7 to 10. One of the interests of the
study is to determine the effect of maternal smoking on the children's
respiratory illness. \cite{laird:84} considered the data collected on
children from Ohio and treated the maternal smoking habit as fixed at
the first visit. The response is binary indicating the presence or
absence of respiratory illness. The maternal smoking habit, in the
preceding year, is recorded as a binary covariate. The mean response
is modeled as a function of Age, Smoking habit and the
interaction. One of the goals of this study was to assess the effect
of maternal smoking on children's respiratory illness. Note that
measurements observed on each child are serially correlated.

We fit the following logistic model to the binary data
\bes
  \log\left\{ \frac{\mu_{it}}{(1 - \mu_{it})} \right\} = \beta_0
  + \beta_1 X_{i1} + \beta_2 X_{i2} + \beta_3 X_{i1} X_{i2}
\ees
for $i = 1, \ldots, 537$ and $t = 0, \ldots, 4$, where $X_{i1}$ and
$X_{i2}$ are the time-independent covariates for the age of the child
and maternal smoking habit, respectively.  The matrix $\bA_i$ is
diagonal with elements $v(\mu_{it}) = \mu_{it} (1 - \mu_{it})$. The
extended score vector $\overline{\bg}_N(\mbt)$ is constructed by
choosing \(s = 2\), \(\bM_1 = \bI\) and \(\bM_2\) as in (\ref{eq:M2}).

The standard errors of $\wmbt$, denoted by $s(\wmbt)$, are computed as
the square root of diagonal elements of $\bJ_N^{-1}(\wmbt)$, where
$\widehat{\bJ}_N(\cdot)$ is defined in (\ref{eq:jhat}). Table
\ref{tab-coef2} presents the estimators  $\wmbt$ and their
corresponding standard errors $s(\wmbt)$, obtained via the QIF library
\citep{loader1:05}.  The \(t\)-ratios suggest that Age is the
significant covariate to include in the model.  The t-ratio for age
has a negative sign indicating that older children are less likely to
have a respiratory illness, whereas mother's smoking habit has a
positive effect on children's respiratory disease, although not
statistically significant. The interaction between the age of the
child and maternal smoking is also not statistically significant.

\begin{table}[htp]
\caption{The parameter estimators and corresponding estimated standard
errors for the Respiratory Health Study under the AR-1 Correlation
Structure.}
\label{tab-coef2}
\vspace{0.25in}
\begin{center}
\begin{tabular}{lrrr} \hline 
Covariates &  $\wmbt$ & $s(\wmbt)$  &  t-ratio \\ \hline 
Intercept & -1.89404 & 0.11903 & -15.91226 \\
Age & -0.12933 & 0.05671 & -2.28062 \\
Smoke & 0.26384 & 0.18962 &  1.39144 \\
Age $\times$ Smoke  & 0.06070 & 0.08791 & 0.69048 \\ \hline 
\end{tabular}
\end{center}
\end{table} 

In order to assess whether the sub-models are adequate, we compute the
QIF under various sub-models with certain parameter restrictions to
perform chi-square tests for comparing different models.

Each row of the Table \ref{tab:laird} represents results for a given
model and the last row provides the full model. For each model, we
compute the parameter estimate \(\wtb\) and report the corresponding
\(Q_N(\wtb)\). The test statistic $T_N$ is obtained via (\ref{eq:TN})
which compares with the full model. The ``df'' column is the
degrees-of-freedom for the test statistic and $P$ is the P-value. From
the table, the models ``Intercept'' and ``Intercept, Smoke'' are
rejected ($P < 0.05$).  The remaining models that include the Age
variable cannot be rejected. 

\begin{table}[htp]
\caption{Testing of hypotheses for the longitudinal data on
children's respiratory disease. The column \(Q_N
\left(\wtb \right)\) is minimum of the QIF obtained under
the submodel, $T_N$ is the value of the test statistic, df is the
degrees of freedom and $P$ is the p-value for the test statistic.}
\vspace{0.5in}
\centerline{
\begin{tabular}{lrccc} \hline
Covariates  & $Q_N\left(\wtb\right)$ & \(T_N\) & df &
$P$\\ \hline  
Intercept            &    11.898  &   7.926  &     3 & 0.048 \\
Intercept, Smoke     &    10.337  &   6.365  &    2  & 0.041 \\
Intercept, Age       &     5.823  &   1.851  &    2  & 0.396 \\
Intercept, Smoke, Age  &   4.449  &   0.477  &    1  & 0.490 \\
Intercept, Smoke, Age, Age $\times$ Smoke   & 3.972 & 0 & - & - \\
\hline  
\end{tabular}} 
\label{tab:laird}
\end{table}

\begin{remark}
The discrepancy between our results and those of QLL is apparently due
to their use of undocumented modifications of the covariance
estimators of \(\cov\{ \overline{\bg}_N(\mbt)\}\). The results in
Table \ref{tab:laird} are based on minimization of \(Q_N(\mbt)\) as
defined in (\ref{eq:qif}).
\end{remark} 

\section{Assessing the Effect of Basis Matrices on the QIF Test}
\label{sec-sim}

In this section we investigate the effect of the choice of basis
matrices, defining the extended score vector
\(\overline{\bg}_N(\mbt)\), on the performance of the QIF test. The
basis matrices are often unknown in advance; hence, it is necessary to
assess the effect of their misspecification. We calculate the size and
power of the QIF test based on $T_N$ through the simulated correlated
data from Bernoulli and Gaussian distributions.

The QIF test statistic is robust to the choice of basis matrices, in
the sense that the null asymptotic chi-squared distribution of $T_N$
is valid whether or not the basis matrices $\bM_j \; (j = 1, \ldots,
s)$ are correctly specified. However, misspecification may adversely
affect the power of the test, since the asymptotic covariance matrix
\(\bJ^{-1}(\mbt_0)\) defined in (\ref{eq:J}) [consequently, the
non-centrality parameter $\delta^2$ in (\ref{eq:nparm})] depends on
the true covariance matrix of $\mbt_0$.

The following is the trade-off for misspecification of the basis
matrices: (1) If too few basis matrices are included in the model,
then the estimator of $\mbt$ may not be efficient; consequently leads
to a loss of power of the test based on $T_N$. (2) If too many basis
matrices are specified in the model, then the dimensionality of the
extended score vector $\overline{\bg}_N(\mbt)$ increases. This can
lead to numerical instability while affecting the power of the QIF
test based on $T_N$ as more components of \(\wC_N(\mbt)\) are being
estimated.

We conducted two simulation experiments each with $N = 50$ subjects
and \(n = 5\) observations per subject to investigate the effect of
the basis matrices on the QIF-based inference.

Let AR-1 refer to representing $\bR^{-1}(\mba) = (\alpha_1 \, \bM_1 +
\alpha_2 \, \bM_2)$ and AR-2 refer to expressing $\bR^{-1}(\mba) = 
(\alpha_1 \, \bM_1 + \alpha_2 \, \bM_2 + \alpha_3 \, \bM_3)$, where
$\bM_1 = \bI$ (the identity matrix of dimension 5),
\[
 \bM_2 = \begin{pmatrix} 0 & 1 & 0 & 0 & 0 \\ 1 & 0 & 1 & 0 & 0 \\ 0 &
  1 & 0 & 1 & 0 \\ 0 & 0 & 1 & 0 & 1 \\ 0 & 0 & 0 & 1 & 0 \\
 \end{pmatrix} \quad \mbox{and} \quad
\hspace{5mm}
 \bM_3 = \begin{pmatrix}
  0 & 0 & 1 & 0 & 0 \\
  0 & 0 & 0 & 1 & 0 \\
  1 & 0 & 0 & 0 & 1 \\
  0 & 1 & 0 & 0 & 0 \\
  0 & 0 & 1 & 0 & 0 \\
  \end{pmatrix}.
\]

\noindent{\em Binary Correlated Responses:} For each subject, binary responses
$Y_{it}$ for $i = 1, \ldots, 50$ and $t = 1, \ldots, 5$ were generated
according to the following two-state Markov chain with the transition
matrix
\bes
  \rho \begin{pmatrix} 1 & 0 \cr 0 & 1 \end{pmatrix}
  + (1 - \rho) \begin{pmatrix} 1 - \mu_i & \mu_i \cr 1 - \mu_i & \mu_i
  \end{pmatrix}, 
\ees
where $\mu_i$ is defined in (\ref{eq:mu}) below. The response vector
\(Y_{it}\) has the stationary distribution \( \begin{pmatrix} 1 -
\mu_i & \mu_i \end{pmatrix} \) and the AR-1 correlation structure
with autocorrelation \(\rho\). We fit the following logistic model to
the binary data
\bea
  \label{eq:mu}
  \log\left\{ \frac{\mu_{i}}{(1 - \mu_{i})} \right\} = \beta_0 +
  \beta_1 \, X_i \quad \mbox{for} \quad i = 1, \ldots, 50,
\eea
where the covariate \(X_i\) is chosen to be equally spaced on the
interval \([-1, 1]\).

For each simulation experiment, the QIF test statistic $T_N$, defined
in (\ref{eq:TN}), for \(H_0\!: \beta_1 = 0\) versus \(H_1\!: \beta_1
\ne 0\) was compared with the critical value \(3.8415\), based on the
95th percentile of the \(\chi_1^2\) distribution. We chose the
following true parameters for the simulation experiment: \(\beta_0 =
\beta_1 = 0\) under \(H_0\) and \(\beta_0 = 0, \beta_1 = 0.5\) under
\(H_1\). 

\noindent{\em Gaussian Correlated Responses:} The same design and
parameter values were used for this simulation experiment; however, we
fit the following model to the continuous data
\bes
  Y_{it} = \beta_0 + \beta_1 X_i + \epsilon_{it} \quad \mbox{for}
  \quad i = 1, \ldots, 50; \; t = 1, \ldots, 5,
\ees
where for each \(i\), \(\epsilon_{it}\) is assumed to be a Gaussian
AR-1 process with variance 1 and correlation \(\rho\).

\begin{table}[htp]
\caption{Achieved significance level and power under three different
assumed correlation structures for the QIF test based on
$T_N$. Results are based on 10,000 replications under the true AR-1
correlation structure with autocorrelation $\rho$.}
\label{tab:simpow}
\vspace{0.5in}
\centerline{
\begin{tabular}{lrccccccc} \hline 
Model & $\rho$ & \multicolumn{3}{c}{Level of Significance} & &
\multicolumn{3}{c}{Power} \\ [2ex]
  & & Identity & AR-1 & AR-2 & & Identity & AR-1 & AR-2 \\ \hline 
         & 0.2 & 0.048 & 0.048 & 0.048  & & 0.473 & 0.463 & 0.447 \\
Logistic & 0.5 & 0.050 & 0.049 & 0.049  & & 0.325 & 0.327 & 0.322 \\ 
         & 0.8 & 0.047 & 0.050 & 0.048  & & 0.226 & 0.228 & 0.227 \\
	 [2ex] 
         & 0.2 & 0.050 & 0.047 & 0.048 & & 0.968 & 0.954 & 0.933 \\
Gaussian & 0.5 & 0.044 & 0.046 & 0.044 & & 0.843 & 0.822 & 0.795 \\
         & 0.8 & 0.052 & 0.050 & 0.051 & & 0.644 & 0.633 & 0.601 \\
	 \hline  
\end{tabular} }
\end{table}

Table \ref{tab:simpow} presents the simulation results under the
logistic and Gaussian models, respectively. Under each scenario, we
achieve significance levels close to the nominal level of 5\%, while
the power decreases as the correlation \(\rho\) increases. However,
for a fixed \(\rho\), there is a minimal difference between the powers
attained under the three different correlation structures. In
particular, there is essentially no power loss when we assumed (albeit
incorrectly) the identity correlation structure.

To investigate further, we estimated the non-centrality parameter 
based on (\ref{eq:nparm}) under different scenarios by finding
\bes
  \widehat{\delta}^{\, 2} = N \, \beta_1 \,
    \left\{ \mbL^T \, \widehat{\bJ}_N^{-1}(\wmbt) \, \mbL
    \right\}^{-1} \, \beta_1,
\ees
averaged over all 10,000 replications, where $\beta_1 = 0.5$ is the
slope parameter under $H_1$, \(\mbL^T = (0 \quad 1)\) and
$\widehat{\bJ}_N^{-1}(\cdot)$ is defined in (\ref{eq:jhat}).

Table \ref{tab:ncppow} presents the $\widehat{\delta}^{\, 2}$ values 
along with the power of the QIF test calculated under the theoretical
asymptotic non-central chi-squared distribution of $T_N$. These
results demonstrate that the model with a misspecified identity
correlation structure yields a slightly smaller $\widehat{\delta}^{\,
2}$ and correspondingly slightly lower power. However, this difference
is not reflected by the finite-sample simulation results presented in
Table
\ref{tab:simpow}.

\begin{table}[htp]
\caption{Estimated non-centrality parameter $\widehat{\delta}^{\, 2}$
values and powers for the QIF test based on \(T_N\) under three
different assumed correlation structures.  Results are based on 10,000
replications under the true AR-1 correlation structure with
autocorrelation $\rho$.}
\label{tab:ncppow}
\vspace{0.5in}
\centerline{
\begin{tabular}{lrccccccc} \hline
Model & $\rho$ & \multicolumn{3}{c}{$\widehat{\delta}^{\, 2}$} & &
\multicolumn{3}{c}{Power} \\ [2ex]
& & Identity & AR-1 & AR-2 & & Identity & AR-1 & AR-2 \\ \hline
         & 0.2 & 4.122 & 4.286 & 4.437 & & 0.528 & 0.544 & 0.558 \\ 
Logistic & 0.5 & 2.452 & 2.609 & 2.678 & & 0.347 & 0.365 & 0.373 \\ 
         & 0.8 & 1.431 & 1.515 & 1.517 & & 0.223 & 0.234 & 0.234 \\
	 [2ex] 
         & 0.2 & 17.982 & 19.263 & 20.478 & & 0.989 & 0.992 & 0.995 \\
Gaussian & 0.5 & 11.067 & 12.189 & 12.961 & & 0.914 & 0.937 & 0.950 \\
         & 0.8 &  6.837 &  7.582 &  8.071 & & 0.744 & 0.786 & 0.811 \\
\hline  
\end{tabular} }
\end{table}

The conclusion from the simulation experiments is that not much is
lost when only the identity correlation structure is assumed. This
does not mean that correlation structure is not important, but rather
that correlation is being adequately modeled through the empirical
covariance matrix \(\wC_N(\wmbt)\) even under the misspecification of
the basis matrices.

\section{Conclusions}
\label{sec-conc}

The QIF is a powerful tool for building regression models for
correlated data. The large-sample properties of the inference
functions are similar to those of the loglikelihood in parametric
statistical inference with test statistics based on the QIF having
asymptotic chi-squared distributions. As shown in Section
\ref{sec-qif}, the covariance matrix for the extended score functions 
defining the QIF employed by QLL can lead to breakdown of the
asymptotic theory.

In this research, we established a unified large-sample theoretical
framework for the QIF. First, we formulated an accurate estimator for
the covariance of the extended score functions
$\overline{\bg}_N(\mbt)$ and second, we derived relevant asymptotic
results necessary for the estimation and testing within the QIF
setting. The key principle underlying our asymptotic treatment is the
quadratic approximation in Theorem \ref{lem:quadap}. The consequences
of this approximation are wide-ranging, providing the necessary
machinery for deriving large-sample theory for the QIF estimators and
test statistics, analogous to standard inferential theory for the
generalized least squares criteria, while leading to a stable and
flexible algorithm for finding the QIF estimators. Our simulation
experiments demonstrate that the QIF test statistic $T_N$ is robust to
the choice of basis matrices, in the sense that the null asymptotic
chi-squared distribution of $T_N$ holds true even under the
misspecification of the basis matrices.

\vspace{0.1in}
{\bf Acknowledgments:} The research of Pilla was supported in part by
the National Science Foundation (NSF) grant DMS 02-39053 and the
Probability and Statistics Program, Office of Naval Research (ONR)
grants N00014-02-1-0316 and N00014-04-1-0481. The research of Loader
was partially supported by the NSF grant DMS 03-06202 and ONR grant
N00014-04-1-0481. The authors thank A.\ Qu for providing the
Respiratory Health Data.

\newpage
\bibliography{qif}
\bibliographystyle{apalike}

\end{document}